\documentclass{article}

\newtheorem{theorem}{Theorem}
\newtheorem{lemma}{Lemma}
\newtheorem{corollary}{Corollary}

\newcommand{\bt}{\begin{theorem}}
\newcommand{\et}{\end{theorem}}
\newcommand{\bl}{\begin{lemma}}
\newcommand{\el}{\end{lemma}}
\newcommand{\bc}{\begin{corollary}}
\newcommand{\ec}{\end{corollary}}
\newcommand{\pf}{{\bf Proof}.\ }
\newcommand{\bq}{\begin{eqnarray*}}
\newcommand{\eq}{\end{eqnarray*}}
\newcommand{\be}{\begin{eqnarray}}
\newcommand{\ee}{\end{eqnarray}}
\newcommand{\beq}{\begin{equation}}
\newcommand{\eeq}{\end{equation}}
\newcommand{\benum}{\begin{enumerate}}
\newcommand{\eenum}{\end{enumerate}}
\newcommand{\ba}{\begin{array}}
\newcommand{\ea}{\end{array}}

\title{Asymptotic density and the asymptotics of partition functions
}
\author{Melvyn B. Nathanson
\thanks{This work was supported in part by grants from the PSC--CUNY Research
Award Program and the NSA Mathematical Sciences Program.}\\
Department of Mathematics\\
Lehman College (CUNY)\\
Bronx, New York 10468\\
e-mail: nathansn@alpha.lehman.cuny.edu}

\date{}

\begin{document}

\maketitle
\begin{abstract}
Let $A$ be a set of positive integers with $\gcd(A) = 1$,
and let $p_A(n)$ be the partition function of $A$.
Let $c_0 = \pi\sqrt{2/3}$.
If $A$ has lower asymptotic density $\alpha$
and upper asymptotic density $\beta$, then
$\liminf \log p_A{n}/c_0\sqrt{n} \geq \sqrt{\alpha}$ and
$\limsup \log p_A(n)/c_0\sqrt{n} \leq \sqrt{\beta}$.
In particular, if $A$ has asymptotic density $\alpha > 0$,
then $\log p_A(n) \sim c_0\sqrt{\alpha n}$.
Conversely, if $\alpha > 0$ and $\log p_A(n) \sim c_0\sqrt{\alpha n}$,
then the set $A$ has asymptotic density $\alpha.$
\end{abstract}

\section{The growth of $p_A(n)$}    \label{para:section:size}
Let $A$ be a nonempty set of positive integers.
The {\em counting function} $A(x)$ of the set $A$
counts the number of positive elements of $A$ that do not exceed $x$.
Then $0 \leq A(x) \leq x$, and so $0 \leq A(x)/x \leq 1$ for all $x$.
The {\em lower asymptotic density} of $A$ is
\[
d_L(A) = \liminf_{x \rightarrow\infty} \frac{A(x)}{x}.
\]
The {\em upper asymptotic density} of $A$ is
\[
d_U(A) = \limsup_{x \rightarrow\infty} \frac{A(x)}{x}.
\]
We have $0 \leq d_L(A) \leq d_U(A) \leq 1$ for every set $A$.
If $d_L(A) = d_U(A)$, then the limit
\[
d(A) = \lim_{x \rightarrow\infty} \frac{A(x)}{x}
\]
exists, and is called the {\em asymptotic density} of the set $A$.

A {\em partition of $n$ with parts in $A$} is a representation of $n$
as a sum of not necessarily distinct elements of $A$, where the number
of summands is unrestricted.
The summands are called the {\em parts} of the partition.
The {\em partition function} $p_A(n)$ counts the number of partitions
of $n$ into parts belonging to the set $A$.
Two partitions that differ only in the order of their parts are
counted as the same partition.
We define $p_A(0) = 1$ and $p_A(-n) = 0$ for $n \geq 1$.

The partition function for the set of all positive integers
is denoted $p(n)$.
Clearly, $0 \leq p_A(n) \leq p(n)$ for every integer $n$
and every set $A$.
A classical result of Hardy and Ramanujan~\cite{hard-rama18}
and Uspensky~\cite{uspe20} states that
\[
\log p(n) \sim c_0\sqrt{n},
\]
where
\[
c_0 = \pi\sqrt{\frac{2}{3}} = 2\sqrt{\frac{\pi^2}{6}}.
\]
Erd\H os~\cite{erdo42b} has given an elementary proof of this result.

Let $\gcd(A)$ denote the greatest common divisor of the elements of $A$.
If $d = \gcd(A) > 1$, consider the set $A' = \{a/d : a\in A\}$.
Then $A'$ is a nonempty set of positive integers such that $\gcd(A') = 1$,
and
\[
p_A(n) = \left\{\ba{ll}
0 & \mbox{if $n\not\equiv 0\pmod {d}$, }\\
p_{A'}\left( n/d\right) & \mbox{if $n\equiv 0\pmod {d}$. }
\ea\right.
\]
Thus, it suffices to consider only  partition functions
for sets $A$ such that $\gcd(A) = 1$.

In this paper we investigate the relationship between
the upper and lower asymptotic densities of a set $A$
and the asymptotic behavior of $\log p_A(n)$.
In particular, we give a complete and elementary proof
of the theorem that, for $\alpha > 0$, the set $A$
has density $\alpha$ if and only if $\log p_A(n) \sim c_0\sqrt{\alpha n}$.
This result was stated, with a sketch of a proof,
in a beautiful paper of Erd\H os~\cite{erdo42b}.

Many other results about the asymptotics of partition functions
can be found in Andrews~\cite[Chapter 6]{andr76} and Odlyzko~\cite{odly95}.

\section{Some lemmas about partition functions}

\bl    \label{para:lemma:monotonic}
Let $A$ be a set of positive integers.
If $p_A(n_0) \geq 1$,
then $p_A(n+n_0) \geq p_A(n)$ for every nonnegative integer $n$.
\el

\pf
The inequality is true for $n=0$, since $p_A(n_0) \geq 1 = p_A(0)$.
We fix one partition $n_0 = a'_1 + \cdots + a'_r$.
Let $n \geq 1$.
To every partition
\[
n = a_1 + \cdots + a_k
\]
we associate the partition
\[
n+n_0 = a_1 + \cdots + a_k + a'_1 + \cdots + a'_r.
\]
This is a one--to--one map from partitions of $n$
to partitions of $n+n_0$, and so $p_A(n) \leq p_A(n+n_0)$.

\bl   \label{para:lemma:monotonic-modm}
Let $A$ be a nonempty set of positive integers,
and let $a_1 \in A$.
For every number $x \geq a_1$ there exists an integer $u$
such that
\[
x-a_1 < u \leq x
\]
and
\[
\max\{p_A(n) : 0 \leq n \leq x\} = p_A(u).
\]
\el

\pf
If $a_1 \in A$, then $p_A(a_1) \geq 1$.
By Lemma~\ref{para:lemma:monotonic},
\[
p_A(n) \leq p_A(n+a_1)
\]
for every nonnegative integer $n$.
Therefore, the partition function $p_A(n)$
is increasing in every congruence class modulo $a_1$.
If $0 \leq r \leq a_1-1$, then
\[
\max\{p_A(n): 0 \leq n \leq x, n\equiv r\pmod{a_1}\} = p_A(u_r)
\]
for some integer $u_r\in (x-a_1,x].$
It follows that
\[
\max\{p_A(n): 0 \leq n \leq x \} = p_A(u),
\]
where
\[
u = \max\{u_0,u_1,\ldots,u_{a_1-1}\} \in (x-a_1,x].
\]
This completes the proof.

\bl                      \label{para:theorem:finiteset}
Let $A$ be a nonempty finite set of relatively prime positive
integers, and let $k$ be the cardinality of $A$.
Let $p_A(n)$ denote the number of partitions of $n$ into parts
belonging to $A$.
Then
\[
p_A(n) = \left(\frac{1}{\prod_{a\in A}a}\right) \frac{n^{k-1}}{(k-1)!}
+ O\left( n^{k-2}\right).
\]
\el

\pf
This is an old result.
The usual proof
(Netto~\cite{nett27}, P\' olya--Szeg\" o~\cite[Problem 27]{poly-szeg25})
is based on the partial fraction
decomposition of the generating function for $p_A(n)$.
There is also an arithmetic proof due to Nathanson~\cite{nath00b}.

\bl      \label{para:lemma:cofinite-1}
Let $n_0$ be a positive integer, and let $A$ be the set of all integers
greater than or equal to $n_0$.
Then $p_A(n)$ is increasing for all positive integers $n$,
and strictly increasing for $n \geq 3n_0+2$.
\el

\pf                                  
If $1 \leq n < n_0$, then $p_A(n) = 0$.

We say that a partition
$a_1 + a_2 + \cdots + a_r$
has a {\em unique largest part} if $a_1 > a_2 \geq \cdots \geq a_r$.
Let $n \geq n_0$.
Then $p_A(n) \geq 1$ since $n \in A$.
To every partition $\pi$ of $n$ we associate a partition of $n+1$
by adding 1 to the largest part of $\pi$.
This is a one--to--one map from the set of all partitions of $n$
and to the set of partitions of $n+1$ with a unique largest part,
and so $p_A(n) \leq p_A(n+1)$ for $n \geq 1$.

Let $n \geq 3n_0+2$.
If $n-n_0$ is even, then $a = (n-n_0)/2 \geq n_0+1$, and $n = 2a+n_0$.
If $n-n_0$ is odd, then $a = (n-n_0-1)/2 \geq n_0+1$, and $n = 2a+(n_0+1)$.
In both cases, $a \in A$.
Therefore, if $n \geq 3n_0+2$, then there exists a partition of $n$
with parts in $A$ and with no unique largest part, and so $p_A(n) < p_A(n+1)$.
This completes the proof.

A set of positive integers is {\em cofinite}
if it contains all  but finitely many positive integers.

\bl \label{para:lemma:cofinite}
Let $A$ be a cofinite set of positive integers.
Then
\[
\log p_A(n) \sim c_0 \sqrt{n}.
\]
\el

\pf
Since $A$ is cofinite,
we can choose an integer $n_0 > 1$ such that $A$ contains
the set $A' = \{n \geq n_0\}$.  Then
\[
p_{A'}(n) \leq p_A(n)\leq p(n).
\]
Since $\log p(n) \sim c_0\sqrt{n}$, it suffices to prove that
$\log p_{A'}(n) \sim c_0\sqrt{n}$.

Let $F = \{n: 1 \leq n \leq n_0-1\}$.
Applying Lemma~\ref{para:theorem:finiteset} with $k = n_0-1$,
we obtain a constant $c \geq 1$
such that $p_F(n) \leq cn^{n_0-2}$ for all positive integers $n$.
Each part of a partition of $n$ must belong  either to $A'$ or to $F$,
and so every partition of $n$ is uniquely of the form $n = n' + (n-n')$,
where $n'$ is a sum of elements of $A'$ and $n-n'$ is a sum of elements of $F$.
By Lemma~\ref{para:lemma:cofinite-1}, the partition function $p_{A'}(n)$
is increasing.
Let $n \geq n_0$.  Then $p_{A'}(n) \geq 1$ and
\bq
p(n) & = & \sum_{n'=0}^n p_{A'}(n')p_F(n-n')  \\
& \leq & cn^{n_0-2}\sum_{n'=0}^n p_{A'}(n')  \\
& \leq & 2cn^{n_0-1} p_{A'}(n).
\eq
Taking logarithms of both sides, we obtain
\bq
\log p(n)
& \leq & \log 2c + (n_0-1)\log n + \log p_{A'}(n)  \\
& \leq & \log 2c + (n_0-1)\log n + \log p(n)
\eq
and so
\bq
\frac{\log p(n)}{c_0\sqrt{n}}
& \leq & \frac{\log 2c + (n_0-1)\log n}{c_0\sqrt{n}}
+ \frac{\log p_{A'}(n)}{c_0\sqrt{n}}  \\
& \leq & \frac{\log 2c + (n_0-1)\log n}{c_0\sqrt{n}}
+ \frac{\log p(n)}{c_0\sqrt{n}}.
\eq
Taking the limit as $n$ goes to infinity,
we have $\log p_{A'}(n) \sim c_0\sqrt{n}$.
This completes the proof.

\section{Abelian and tauberian theorems}
In this section we derive two results in analysis that will be used
in the proof of Theorem~\ref{para:theorem:inverse}.
To every sequence $B = \{b_n\}_{n=0}^{\infty}$ of real numbers
we can associate the power series $f(x) = \sum_{n=0}^{\infty}b_nx^n$.
We shall assume that the power series converges for $|x| < 1$.
We think of the function $f(x)$ as a kind of average over the sequence
$B$.
Roughly speaking, an {\em abelian theorem}\index{abelian theorem}
asserts that if the sequence $B$ has some property, then the function
$f(x)$ has some related property.
Conversely, a {\em tauberian theorem}\index{tauberian theorem}
asserts that if the function $f(x)$ has some property,
then the sequence $B$ has a related property.

The following theorem is abelian.

\bt       \label{para:theorem:abelian}
Let $B = \{ b_n\}_{n=0}^{\infty}$ be a sequence
of nonnegative real numbers
such that the power series
$f(x) = \sum_{n=0}^{\infty} b_nx^n$
converges for $|x| < 1$.
If
\beq                        \label{para:liminfalpha}
\liminf_{n \rightarrow\infty} \frac{\log b_n}{\sqrt{n}} \geq 2\sqrt{\alpha},
\eeq
then
\beq                        \label{para:liminf}
\liminf_{x \rightarrow 1^-} (1-x)\log f(x) \geq \alpha.
\eeq
If
\beq                        \label{para:limsupbeta}
\limsup_{n \rightarrow\infty} \frac{\log b_n}{\sqrt{n}} \leq 2\sqrt{\beta},
\eeq
then
\beq                        \label{para:limsup}
\limsup_{x \rightarrow 1^-} (1-x)\log f(x) \leq \beta.
\eeq
In particular, if $\alpha > 0$ and
\beq                        \label{para:2-x}
\log b_n \sim 2\sqrt{\alpha n},
\eeq
then
\beq      \label{para:1-x}
\log f(x) \sim \frac{\alpha}{1-x}.
\eeq
\et

\pf
Let $0 < \varepsilon < 1.$
Inequality~(\ref{para:liminfalpha}) implies that there exists
a positive integer $N_0 = N_0(\varepsilon)$ such that
\[
b_n > e^{2(1-\varepsilon)\sqrt{\alpha n}} \quad\mbox{for all $n \geq N_0$}.
\]
For $0 < x < 1$, we let $x = e^{-t}$,
where $t = t(x) = -\log x > 0$, and $t$ decreases to 0
as $x$ increases to 1.

If $n \geq N_0$, then
\[
b_nx^n > e^{2(1-\varepsilon)\sqrt{\alpha n}}e^{-tn}
= e^{2(1-\varepsilon)\sqrt{\alpha n}-tn}.
\]
Completing the square in the exponent, we obtain
\[
2(1-\varepsilon)\sqrt{\alpha n}-tn
= \frac{(1-\varepsilon)^2\alpha}{t}
  - t\left( \sqrt{n}-\frac{(1-\varepsilon)\sqrt{\alpha}}{t} \right)^2,
\]
and so
\[
b_nx^n > e^{\frac{(1-\varepsilon)^2\alpha}{t}}
e^{ - t\left( \sqrt{n}-\frac{(1-\varepsilon)\sqrt{\alpha}}{t} \right)^2}.
\]
Choose $t_0 > 0$ such that
\[
\frac{(1-\varepsilon)^2 \alpha}{t_0^2}  > N_0 + 1,
\]
and let $x_0 = e^{-t_0} < 1$.
If $x_0 < x < 1$ and $x = e^{-t}$, then $0 < t < t_0$.
Let
\[
n_x = \left[ \frac{(1-\varepsilon)^2 \alpha}{t^2} \right].
\]
Then
\[
N_0 < \frac{(1-\varepsilon)^2 \alpha}{t^2} - 1
< n_x \leq \frac{(1-\varepsilon)^2 \alpha}{t^2}
\]
and
\[
\frac{(1-\varepsilon)\sqrt{\alpha}}{t} - 1
< \sqrt{\frac{(1-\varepsilon)^2 \alpha}{t^2} - 1 }
< \sqrt{n_x} \leq \frac{(1-\varepsilon)\sqrt{\alpha}}{t}.
\]
It follows that
\[
\left( \sqrt{n_x} - \frac{(1-\varepsilon)\sqrt{\alpha}}{t} \right)^2 < 1,
\]
and so
\[
b_{n_x}x^{n_x} > e^{\frac{(1-\varepsilon)^2\alpha}{t}}e^{-t}
= e^{\frac{(1-\varepsilon)^2\alpha}{t} -t}.
\]
Since $b_nx^n \geq 0$ for all $n \geq 0$, we have
\[
f(x) = \sum_{n=0}^{\infty} b_nx^n
\geq b_{n_x}x^{n_x} > e^{\frac{(1-\varepsilon)^2\alpha}{t} -t}.
\]
Therefore,
\[
\log f(x) > \frac{(1-\varepsilon)^2\alpha}{t} -t
\]
and
\[
t\log f(x) > (1-\varepsilon)^2\alpha - t^2.
\]                                           
Since
\[
t = - \log x \sim 1-x \quad\mbox{as $x \rightarrow 1^-$,}
\]
it follows that
\bq
\liminf_{x\rightarrow 1^-}(1-x)\log f(x)
& = & \liminf_{x\rightarrow 1^-} t\log f(x) \\
& \geq & \liminf_{t\rightarrow 0^+}
\left( (1-\varepsilon)^2\alpha - t^2 \right)  \\
& = & (1-\varepsilon)^2\alpha .
\eq
This inequality is true for every $\varepsilon > 0,$
and so
\[
\liminf_{x\rightarrow 1^-}(1-x)\log f(x) \geq \alpha .
\]
This proves~(\ref{para:liminf}).

If~(\ref{para:limsupbeta}) holds, then there exists a positive integer
$N_0 = N_0(\varepsilon)$ such that
\[
b_n < e^{2(1+\varepsilon)\sqrt{\beta n}} \quad\mbox{for all $n \geq N_0$}
\]
Let $x = e^{-t}$.  Then
\bq
f(x) & = & \sum_{n=0}^{\infty} b_nx^n  \\
& < & \sum_{n=0}^{N_0-1} b_nx^n
+ e^{\frac{(1+\varepsilon)^2\beta}{t}} \sum_{n=N_0}^{\infty}
e^{ - t\left( \sqrt{n}-\frac{(1+\varepsilon)\sqrt{\beta}}{t} \right)^2}  \\
& = & c_1(\varepsilon)
+ e^{\frac{(1+\varepsilon)^2\beta}{t}} \sum_{n=N_0}^{\infty}
e^{ - t\left( \sqrt{n}-\frac{(1+\varepsilon)\sqrt{\beta}}{t} \right)^2},
\eq
where
\[
0 \leq \sum_{n=0}^{N_0-1} b_n x^n
\leq \sum_{n=0}^{N_0-1} b_n = c_1(\varepsilon).
\]
If
\[
n > \left[\frac{16\beta}{t^2}\right] = N_1(t) = N_1,
\]
then
\[
\sqrt{n} > \frac{4\sqrt{\beta}}{t} > \frac{2(1+\varepsilon)\sqrt{\beta}}{t}
\]
and
\[
\sqrt{n} - \frac{(1+\varepsilon)\sqrt{\beta}}{t} > \frac{\sqrt{n}}{2}.
\]
It follows that
\[
e^{ - t\left( \sqrt{n}-\frac{(1+\varepsilon)\sqrt{\beta}}{t} \right)^2}
< e^{ - t\left( \frac{\sqrt{n}}{2} \right)^2}
= e^{-\frac{tn}{4}},
\]
and so
\bq
\sum_{n=N_1+1}^{\infty}
e^{ - t\left( \sqrt{n}-\frac{(1+\varepsilon)\sqrt{\beta}}{t} \right)^2}
& < & \sum_{n=N_1+1}^{\infty} e^{-\frac{tn}{4}}  \\
& = & \frac{e^{ - t(N_1+1)/4 }}{1 - e^{ - t/4 } }  \\
& < & \frac{8e^{ - 4\beta/t }}{t},
\eq
since $1-t/4 < e^{-t/4} < 1- t/8$ for $0 < t < 1$.
Moreover,
\[
\sum_{n= N_0}^{N_1}
e^{ - t\left( \sqrt{n}-\frac{(1+\varepsilon)\sqrt{\beta}}{t} \right)^2}
< N_1 \leq \frac{16\beta}{t^2}.
\]
Consequently,
\bq
f(x) & \leq & c_1(\varepsilon)
+ e^{\frac{(1+\varepsilon)^2\beta}{t}}
\left( \frac{16\sqrt{\beta}}{t^2} + \frac{8e^{ - 4\beta/t }}{t} \right)  \\
& \leq & \frac{c_2(\varepsilon)e^{\frac{(1+\varepsilon)^2\beta}{t}}}{t^2}.
\eq
Therefore,
\[
\log f(x) \leq \frac{(1+\varepsilon)^2\beta}{t}
+ \log \frac{c_2(\varepsilon)}{t^2},
\]
and so
\[
t\log f(x) \leq (1+\varepsilon)^2\beta + t\log \frac{c_2(\varepsilon)}{t^2}.
\]
Then
\[
\limsup_{x\rightarrow 1^-} (1-x)\log f(x)
= \limsup_{t\rightarrow 0^+} t\log f(x) \leq (1+\varepsilon)^2\beta.
\]
This inequality is true for every $\varepsilon > 0$, and so
\[
\limsup_{x\rightarrow 1^-} (1-x)\log f(x) \leq \beta.
\]
This proves~(\ref{para:limsup}).

If~(\ref{para:2-x}) holds, that is, if
\[
\lim_{n \rightarrow\infty} \frac{\log b_n}{2\sqrt{n}} = \sqrt{\alpha} > 0,
\]
then~(\ref{para:liminfalpha}) and~(\ref{para:limsupbeta})
hold with $\alpha = \beta$.
These inequalities imply~(\ref{para:liminf}) and~(\ref{para:limsup}),
and so
\[
\lim_{x\rightarrow 1^-} (1-x)\log f(x) = \alpha,
\]
or, equivalently,
\[
\log f(x) \sim \frac{\alpha}{1-x}.
\]
This completes the proof.

The statement that~(\ref{para:2-x}) implies~(\ref{para:1-x})
appears in Erd\H os~\cite{erdo42b}.

The following tauberian theorem generalizes a
well--known result of Hardy and Littlewood~\cite{hard-litt14}.

\bt       \label{para:theorem:tauberian}
Let $B = \{b_n\}_{n=0}^{\infty}$ be a sequence of nonnegative real numbers
such that the power series
\[
f(x) = \sum_{n=0}^{\infty} b_nx^n
\]
converges for $|x| < 1$.
Let
\[
S_B(n) = \sum_{k=0}^n b_k.
\]
Let $c > 0.$
If
\beq     \label{para:tauberian1}
\limsup_{x \rightarrow 1^{-}} (1-x)f(x) \leq c,
\eeq
then
\beq     \label{para:tauberian2}
\limsup_{n \rightarrow \infty} \frac{S_B(n)}{n} \leq c.
\eeq
If
\beq     \label{para:tauberian3}
\liminf_{x \rightarrow 1^{-}} (1-x)f(x) \geq c,
\eeq
then
\beq     \label{para:tauberian4}
\liminf_{n \rightarrow \infty} \frac{S_B(n)}{n} \geq c.
\eeq
In particular, if
\beq     \label{para:tauberian5}
f(x) \sim \frac{c}{1-x} \qquad\mbox{as $x \rightarrow 1^-$,}
\eeq
then
\beq     \label{para:tauberian6}
S_B(n) \sim cn.
\eeq
\et

\pf
The Hardy--Littlewood theorem states
that~(\ref{para:tauberian5}) implies~(\ref{para:tauberian6}).
The proofs that~(\ref{para:tauberian1}) implies~(\ref{para:tauberian2})
and that~(\ref{para:tauberian3}) implies~(\ref{para:tauberian4})
require only a simple modification of Karamata's method,
as presented in Titchmarsh~\cite[Chapter 7]{titc39}.

\section{Direct and inverse theorems for $p_A(n)$}
A {\em direct theorem} uses information about the sequence $A$
to deduce properties of the partition function $p_A(n)$.
An {\em inverse theorem} uses information about the partition
function $p_A(n)$ to deduce properties of the sequence $A$.
We begin with a direct theorem.

\bt                    \label{para:theorem:direct}
Let $A$ be an infinite set of positive integers with $\gcd(A) = 1$.
If $d_L(A) \geq \alpha$, then
\[
\liminf_{n \rightarrow\infty}\frac{\log p_A(n)}{c_0\sqrt{n}} \geq \sqrt{\alpha}.
\]
If $d_U(A) \leq \beta$, then
\[
\limsup_{n \rightarrow\infty}\frac{\log p_A(n)}{c_0\sqrt{n}} \leq \sqrt{\beta}.
\]
\et

\pf
Let $A = \{a_k\}_{k=1}^{\infty}$, where $a_1 < a_2 < \cdots$.
Since $\gcd(A) = 1$,
there is an integer $\ell_0$ such that
$\gcd\{a_k: 1 \leq k \leq \ell_0 -1 \} = 1$.
Let $\varepsilon > 0$.
If $d_U(A) \leq \beta$, there exists an integer
$k_0 = k_0(\varepsilon) \geq \ell_0$ such that,
for all $k \geq k_0$,
\[
\frac{k}{a_k} = \frac{A(a_k)}{a_k} < \beta + \varepsilon,
\]
and so
\[
k < (\beta + \varepsilon)a_k.
\]
Let $A' = \{a_k\in A : k \geq k_0 \}$
and $F = A \setminus A' = \{a_k\in A : 1 \leq k \leq k_0 -1 \}$.

Let $n $ and $n'$ be positive integers, $n' \leq n$,
and let
\[
n' = a_{k_1} + a_{k_2} + \cdots + a_{k_r}
\]
be a partition of $n'$ with parts in $A'$.
Then $k_i \geq k_0$ for all $i = 1, \ldots, r$.
To this partition of $n'$ we associate the partition
\[
m = k_1 + k_2 + \cdots + k_r.
\]
Since $k_i < (\beta + \varepsilon)a_{k_i}$ for $i = 1,\ldots, r$,
we have
\bq
m & < & (\beta + \varepsilon)a_{k_1} + (\beta + \varepsilon)a_{k_2}
+ \cdots + (\beta + \varepsilon) a_{k_r}  \\
& = & (\beta + \varepsilon)n'  \\
& \leq & (\beta + \varepsilon)n.
\eq
This is a one--to--one mapping from partitions of $n'$
with parts in $A'$ to partitions of integers less than
$(\beta+\varepsilon)n$, and so
\bq
p_{A'}(n')
& \leq & \sum_{m < (\beta+\varepsilon)n } p(m)  \\
& \leq & (\beta+\varepsilon)n \max\{p(m) : m < (\beta+\varepsilon)n \}\\
& \leq & (\beta+\varepsilon)n p( [(\beta+\varepsilon)n]) \\
&  <   & 2n p( [(\beta+\varepsilon)n]),
\eq
since the unrestricted partition function $p(n)$ is strictly increasing.

We have $A = A' \cup F$, where $A' \cap F = \emptyset$.
The set $F$ is a nonempty finite set of integers
of cardinality $k_0 -1$, and $\gcd(F) = 1$ since $k_0 \geq \ell_0$.
By Theorem~\ref{para:theorem:finiteset}, there exists a constant $c$
such that
\[
p_F(n) \leq cn^{k_0-2}
\]
for every positive integer $n$.
Every partition of $n$ with parts in $A$ can be decomposed uniquely into
a partition of $n'$ with parts in $A'$ and a partition of $n-n'$
with parts in $F$, for some nonnegative integer $n' \leq n$.
Then
\bq
p_A(n)
& = & \sum_{n' = 0}^n p_{A'}(n')p_F(n-n')  \\
& \leq & cn^{k_0-2}\sum_{n' = 0}^n p_{A'}(n')  \\
& \leq & cn^{k_0-2}(n+1)\max\{p_{A'}(n'): n' = 0,1,\ldots, n\}  \\
& \leq & 2cn^{k_0-1}\max\{p_{A'}(n'): n' = 0,1,\ldots, n\}  \\
& < & 2cn^{k_0-1}2n p( [(\beta+\varepsilon)n])  \\
& = & 4cn^{k_0} p( [(\beta+\varepsilon)n]) .
\eq
Since $\log p(n) \sim c_0\sqrt{n}$, it follows that for every $\varepsilon>0$
there exists an integer $n_0(\varepsilon)$ such that
\[
\log p(n) < (1+\varepsilon)c_0\sqrt{n}
\]
for $n \geq n_0(\varepsilon).$
Therefore,
\bq
\log p_{A}(n)
& \leq & \log 4c + k_0\log n + \log p( [(\beta+\varepsilon)n]) \\
& < & \log 4c+k_0\log n+(1+\varepsilon)c_0 \sqrt{(\beta+\varepsilon)n}
\eq
for $n \geq (n_0(\varepsilon)+1)/(\beta + \varepsilon).$
Dividing by $c_0\sqrt{n}$, we obtain
\[
\frac{\log p_A(n)}{c_0\sqrt{n}}
\leq \frac{\log 4c+k_0\log n}{c_0\sqrt{n}}
+ (1+\varepsilon)\sqrt{\beta +\varepsilon},
\]
and so
\[
\limsup_{n\rightarrow \infty} \frac{\log p_A(n)}{c_0\sqrt{n}}
\leq (1+\varepsilon)\sqrt{\beta+\varepsilon}.
\]
Since this inequality is true for all $\varepsilon > 0$, we obtain
\[
\limsup_{n\rightarrow \infty} \frac{\log p_A(n)}{c_0\sqrt{n}} \leq \sqrt{\beta}.
\]

Next we prove that if $d_L(A) \geq \alpha$, then
\[
\liminf_{n\rightarrow\infty} \frac{\log p_A(n)}{c_0\sqrt{n}}
\geq \sqrt{\alpha}.
\]
This inequality is trivial if $\alpha = 0$,
since $\log p_A(n)/c_0\sqrt{n} \geq 0$ for all sufficiently large $n$.

Let $\alpha > 0$ and
\[
0 < \varepsilon < \alpha.
\]
There exists an integer
$k_0 = k_0(\varepsilon)$ such that,
for all $k \geq k_0$,
\[
\frac{k}{a_k} = \frac{A(a_k)}{a_k} > \alpha - \varepsilon,
\]
and so
\[
a_k < \frac{k}{\alpha - \varepsilon}.
\]
Since $\gcd(A) = 1$,
every sufficiently large integer can be written as a sum of elements of $A$,
and so there exists an integer $N_0$ such that
$p_A(n) \geq 1$ for all $n \geq N_0$.
Let $p'(n)$ denote the number of partitions of $n$
into parts $k \geq k_0$.
To every partition
\[
n = k_1 + \cdots + k_r \qquad \mbox{with $k_1 \geq \cdots \geq k_r \geq k_0$,}
\]
we associate the partition
\[
m = a_{k_1} + \cdots + a_{k_r}.
\]
Then
\[
m < \frac{k_1}{\alpha - \varepsilon} + \cdots
+ \frac{k_1}{\alpha - \varepsilon}
= \frac{n}{\alpha - \varepsilon}.
\]
This is a one--to--one mapping from partitions of $n$
with parts greater than or equal to $k_0$
to partitions of integers $m$ less than
$n/(\alpha - \varepsilon)$, and so
\bq
p'(n) & \leq &
\sum_{m <\frac{n}{\alpha -\varepsilon}} p_A(m)\\
& \leq & \frac{n}{\alpha - \varepsilon}
\max\left\{ p_A(m): m <\frac{n}{\alpha -\varepsilon} \right\} \\
& < & \frac{n}{\alpha - \varepsilon}p_A(u_n),
\eq
where, by Lemma~\ref{para:lemma:monotonic-modm}
(since $a_1 \in A$),
the integer $u_n$ belongs to the bounded interval
\[
\frac{n}{\alpha -\varepsilon} - a_1 < u_n \leq \frac{n}{\alpha -\varepsilon}.
\]
The sequence $\{u_n\}_{n=1}^{\infty}$
is not necessarily increasing, but
\[
\lim_{n\rightarrow\infty} u_n = \infty.
\]
Let $d$ be the unique positive integer such that
\[
0 < (\alpha - \varepsilon)a_1 \leq d < (\alpha - \varepsilon)a_1 + 1.
\]
For every $i,j \geq 1$,
\bq
u_{(i+j)d} - u_{id}
& > & \left( \frac{(i+j)d}{\alpha - \varepsilon} -a_1\right)
 -  \frac{id}{\alpha - \varepsilon}                         \\
& = & \frac{jd}{\alpha - \varepsilon} - a_1                   \\
& \geq & (j-1)a_1.
\eq
It follows that $u_{(i+1)d} > u_{id}$,
and so the sequence $\{u_{id}\}_{i=1}^{\infty}$
is strictly increasing.
Similarly,
\bq
u_{(i+j)d} - u_{id}
& < & \frac{(i+j)d}{\alpha - \varepsilon}
- \left(\frac{id}{\alpha - \varepsilon} - a_1\right)  \\
& = & \frac{jd}{\alpha - \varepsilon} + a_1               \\
& < & (j+1)a_1 + \frac{j}{\alpha - \varepsilon}.
\eq
Let $j_0$ be the unique integer such that
\[
\frac{N_0}{a_1}+1 \leq j_0 < \frac{N_0}{a_1}+2.
\]
Then
\[
u_{i d} - u_{(i - j_0)d} > (j_0-1)a_1 \geq N_0
\]
for all $i \geq j_0$.

For every integer $n \geq j_0 d$ there exists a unique integer
$\ell \geq j_0$ such that
\[
u_{\ell d} \leq n < u_{(\ell +1)d}.
\]
Then
\[
n - u_{(\ell - j_0)d} < u_{(\ell +1) d} - u_{(\ell - j_0)d}
<  (j_0+2)d + \frac{j_0+1}{\alpha - \varepsilon}
\]
and
\[
n - u_{(\ell - j_0)d} \geq u_{\ell d} - u_{(\ell - j_0)d} > N_0.
\]
Since
\[
p_A(n - u_{(\ell - j_0)d}) \geq 1,
\]
Lemma~\ref{para:lemma:monotonic} implies that
\[
p_A(n) \geq p_A(u_{(\ell - j_0)d})
> \left( \frac{\alpha - \varepsilon}{(\ell - j_0)d} \right) p'((\ell - j_0)d).
\]
Since
\[
n < u_{(\ell + 1)d} \leq \frac{(\ell +1)d}{\alpha - \varepsilon},
\]
it follows that
\[
(\ell - j_0)d > (\alpha - \varepsilon)n - (j_0+1)d.
\]
Since $p'(n)$ is the partition function of a cofinite subset of the
positive integers, Lemma~\ref{para:lemma:cofinite}
implies that for $n$ sufficiently large,
\bq
\log p_A(n)
& > & \log p'((\ell - j_0)d) + \log(\alpha - \varepsilon) - \log (\ell - j_0)d  \\
& > & (1-\varepsilon)c_0\sqrt{(\ell - j_0)d} + \log(\alpha - \varepsilon)  - \log (\ell - j_0)d  \\
& > & (1-\varepsilon) c_0\sqrt{(\alpha - \varepsilon)n - (j_0+1)d}
+ \log(\alpha - \varepsilon)- \log (\ell - j_0)d.
\eq
Dividing by $c_0\sqrt{n}$, we obtain
\[
\liminf_{n\rightarrow\infty} \frac{\log p_A(n)}{c_0\sqrt{n}}
\geq (1-\varepsilon) \sqrt{\alpha - \varepsilon}.
\]
This inequality holds for  $0 < \varepsilon < \alpha$,
and so
\[
\liminf_{n\rightarrow\infty} \frac{\log p_A(n)}{c_0\sqrt{n}}
\geq \sqrt{\alpha}.
\]
This completes the proof.

\bt   \label{para:theorem:direct2}
Let $A$ be a set of positive integers with $\gcd(A) = 1$.
If $d(A) = \alpha > 0$, then $\log p_A(n) \sim c_0\sqrt{\alpha n}$.
\et

\pf
This follows from Theorem~\ref{para:theorem:direct}
with $\alpha = \beta$.

\bt        \label{para:theorem:arithpro}
Let $a_1,\ldots, a_{\ell}, m$ be integers such that
\[
1 \leq a_1 < \cdots < a_{\ell} \leq m
\]
and
\[
(a_1,\ldots, a_{\ell},m) = 1.
\]
Let $A$ be the set of all positive integers $a$ such that
$a\equiv a_i\pmod{m}$ for some $i = 1,\ldots, \ell$.
Then
\[
\log p_A(n) \sim c_0 \sqrt{\frac{\ell n}{m}}.
\]
\et

\pf
The set $A$ satisfies $\gcd(A) = 1$ and $d(A) = \ell/m$,
and so the result follows from Theorem~\ref{para:theorem:direct2}
with $\alpha = \ell/m$.
Using Erd\H os's elementary method, Nathanson~\cite{nath99g}
has also given a direct proof of Theorem~\ref{para:theorem:arithpro}.

\bt
Let $A$ be a set of positive integers with $\gcd(A) = 1$.
If $d(A) = 0$, then $\log p_A(n) = o(\sqrt{n}).$
\et

\pf
If $A$ is infinite, this follows from
Theorem~\ref{para:theorem:direct} with $\beta = 0$.
If $A$ is finite, this follows from Lemma~\ref{para:theorem:finiteset}.

The next result is an inverse theorem; it shows how the growth of the
partition function $p_A(n)$ determines the asymptotic
density of the sequence $A$.

\bt                      \label{para:theorem:inverse}
Let $A$ be an infinite set of positive integers with $\gcd(A) = 1$.
If $\alpha > 0$ and
\beq       \label{para:asymptoticformula}
\log p_A(n) \sim c_0\sqrt{\alpha n} = 2\sqrt{\frac{\pi^2\alpha n}{6}},
\eeq
then $A$ has asymptotic density $\alpha$.
\et

\pf
The generating function
\[
f(x) = \sum_{n=0}^{\infty} p_A(n)x^n
= \prod_{a\in A} (1-x^a)^{-1}
\]
converges for $|x|<1$, and
\bq
\log f(x)
& = & -\sum_{a\in A} \log(1-x^a)  \\\
& = & \sum_{a\in A} \sum_{k=1}^{\infty} \frac{x^{ak}}{k} \\
& = & \sum_{\ell=1}^{\infty} b_{\ell} x^{\ell},
\eq
where
\[
b_{\ell} = \sum_{a\in A\atop {\ell} = ak} \frac{1}{k} \geq 0.
\]
Let
\[
S_B(x) = \sum_{\ell \leq x} b_{\ell}.
\]
Then $S_B(x) \geq 0$ for all $x$,
and $S_B(x) = 0$ if $x < 1$.
We have
\bq
S_B(n) & = & \sum_{\ell =1}^n b_{\ell}
= \sum_{\ell =1}^n \sum_{a\in A\atop {\ell} = ak} \frac{1}{k}  \\
& = & \sum_{k =1}^n \frac{1}{k}\sum_{a\in A\atop {a \leq n/k}} 1
= \sum_{k =1}^n \frac{1}{k} A\left(\frac{n}{k}\right).
\eq
By M\" obius inversion, we have
\[
A(n) = \sum_{k =1}^n \frac{\mu(k)}{k} S_B\left(\frac{n}{k}\right).
\]
By Theorem~\ref{para:theorem:abelian},
the asymptotic formula~(\ref{para:asymptoticformula}) implies that
\[
(1-x)\log f(x) \sim \frac{\pi^2\alpha}{6}
\quad\mbox{ as $x\rightarrow 1^-$.}
\]
Theorem~\ref{para:theorem:tauberian} implies that
\[
S_B(n) \sim \frac{\pi^2\alpha n}{6}.
\]
We define the function $r(x)$ by
\[
\frac{S_B(x)}{x} = \frac{\pi^2\alpha}{6} + r(x).
\]
Then $r(x) = o(x)$.
For every $\varepsilon > 0$ there exists an integer
$n_0 = n_0(\varepsilon) > e^2$ such that
\[
|r(x)| < \varepsilon
\]
for all $x \geq n_0.$
If $k > n/n_0$, then $n/k < n_0$ and $0 \leq S_B(n/k) \leq S_B(n_0)$.
Therefore,
\bq
A(n)
& = & \sum_{k =1}^{n} \frac{\mu(k)}{k} S_B\left(\frac{n}{k}\right) \\
& = & \sum_{1 \leq k \leq n/n_0} \frac{\mu(k)}{k} S_B\left(\frac{n}{k}\right)
+ \sum_{n/n_0 < k \leq n} \frac{\mu(k)}{k} S_B\left(\frac{n}{k}\right) \\
& = & \frac{\pi^2\alpha n}{6}\sum_{1 \leq k \leq n/n_0} \frac{\mu(k)}{k^2}
+ n\sum_{1 \leq k \leq n/n_0} \frac{\mu(k)}{k^2} r\left(\frac{n}{k}\right)  \\
& & + \sum_{n/n_0 < k \leq n} \frac{\mu(k)}{k} S_B\left(\frac{n}{k}\right).
\eq
We evaluate these three terms separately.
Since
\[
\sum_{1 \leq k \leq n/n_0} \frac{\mu(k)}{k^2}
= \frac{6}{\pi^2} - \sum_{k > n/n_0} \frac{\mu(k)}{k^2}
= \frac{6}{\pi^2} + O\left(\frac{n_0}{n}\right),
\]
it follows that
\[
\frac{\pi^2\alpha n}{6}\sum_{1 \leq k \leq n/n_0} \frac{\mu(k)}{k^2}
= \alpha n + O\left(1\right).
\]
The second term satisfies
\[
\left|
n\sum_{1 \leq k \leq n/n_0} \frac{\mu(k)}{k^2} r\left(\frac{n}{k}\right)\right|
\leq  \varepsilon n\sum_{1 \leq k \leq n/n_0} \frac{1}{k^2}
= O(\varepsilon n).
\]
The last term is bounded independent of $n$, since
\[
\left|\sum_{n/n_0 < k\leq n}\frac{\mu(k)}{k} S_B\left(\frac{n}{k}\right)\right|
\leq S_B(n_0) \sum_{n/n_0 < k\leq n}\frac{1}{k}
\leq 2 S_B(n_0) \log n_0
= O(1).
\]
Therefore,
\[
A(n) = \alpha n + O(\varepsilon n) + O(1),
\]
and so $d(A) = \alpha$.
This completes the proof.

\bt        \label{para:theorem:erdos}
Let $A$ be a set of positive integers with $\gcd(A) = 1$,
and let $\alpha > 0$.
Then $d(A) = \alpha$ if and only if
\[
\log p_A(n) \sim c_0\sqrt{\alpha n}.
\]
\et

\pf
This follows immediately from Theorem~\ref{para:theorem:direct}
and Theorem~\ref{para:theorem:inverse}.

{\bf Remark}.  Let $A$ be an infinite set of positive integers
with $\gcd(A) = 1$.
Let $\alpha$ and $\beta$ be nonnegative real numbers such that
\[
\liminf_{n \rightarrow\infty} \frac{\log p_A(n)}{c_0\sqrt{n}} \geq \sqrt{\alpha}
\]
and
\[
\limsup_{n \rightarrow\infty} \frac{\log p_A(n)}{c_0\sqrt{n}} \leq \sqrt{\beta}.
\]
Does it follow that $d_L(A) \geq \alpha$ and $d_U (A) \leq \beta$?
This would imply that
$d_L(A) = \alpha$ if and only if
$\liminf_{n \rightarrow\infty} \log p_A(n)/c_0\sqrt{n} = \sqrt{\alpha}$,
and $d_U(A) = \beta$ if and only if
$\limsup_{n \rightarrow\infty} \log p_A(n)/c_0\sqrt{n} = \sqrt{\beta}$.

\end{document}